\documentclass[a4paper,11pt]{article}
\usepackage{amssymb,amsmath,amsthm}
\usepackage{texdraw}
\usepackage{epic,eepic}

\newtheorem{thm}{Theorem}[section]
\newtheorem{prop}[thm]{Proposition}

\theoremstyle{remark}
\newtheorem{rem}[thm]{Remark}

\begin{document}

\title{Yang-Baxter maps and multi-field integrable lattice equations}
\author{Vassilios G. Papageorgiou{\footnote{\tt vassilis@math.upatras.gr}} $\phantom{.}$
and Anastasios G. Tongas{\footnote{\tt tasos@math.upatras.gr}} \\ {\em
Department of Mathematics, University of Patras, 265 00 Patras, Greece} }

\maketitle

\begin{abstract}
A variety of Yang-Baxter maps
are obtained from integrable multi-field equations on quad-graphs. 
A systematic framework for investigating
this connection relies on the symmetry groups of the equations.
The method is applied to lattice equations introduced by 
Adler and Yamilov and which are related to the nonlinear superposition 
formulae for the B\"acklund transformations of the nonlinear Schr\"odinger 
system and specific ferromagnetic models.
\end{abstract}

\section{Introduction}

The study of set-theoretic solutions of the quantum Yang-Baxter (YB) equation
was originally suggested by Drinfeld in \cite{drin} and has attracted 
the interest of many researchers. 
Certain examples of such solutions had already appeared in the 
relevant literature by Sklyanin \cite{Skl}. 
Set-theoretic solutions of the quantum YB equation, termed as classical solutions 
of the quantum YB equation, were investigated by 
Weinstein and Xu \cite{WX} from a geometric point of view using the theory of 
Poisson Lie groups and the associated symplectic groupoids. 
The algebraic theory of the subject was developed by Etingof, Schedler and 
Soloviev \cite{ESS}. 
More recently, the dynamical theory of the set-theoretic solutions of the quantum 
YB equation was developed by Veselov 
in \cite{ves1} (see also the review article \cite{ves2}), 
where the notion of the associated transfer maps was introduced 
as a dynamical analogue of the transfer matrices of the theory of solvable 
models in statistical mechanics. 
Veselov suggested the short term YB maps when referring to set-theoretic 
(or classical) solutions of the YB equation and which we adopt in the present work.

Let $X$ be a set. A map $R:X \times X \rightarrow X \times X$, is called YB 
map if it satisfies the relation
\begin{equation}
R^{(2,3)}\, R^{(1,3)}\, R^{(1,2)}\, = R^{(1,2)}\, R^{(1,3)}\, R^{(2,3)}\, , 
\label{eq:YBrel}
\end{equation}
regarded as an equality of composite maps from $X \times X \times X$ into itself.
Following the usual notation, $R^{(1,3)}$ is meant as the identity in the second 
factor of $ X \times X\times X$ and as $R$ in the first 
and third factors, and analogously for $R^{(1,2)}$, $R^{(2,3)}$ (see e.g. \cite{ves2}).
In various interesting examples of YB maps, such as maps arising from geometric 
crystals \cite{eting2}, the set $X$ has the structure of an algebraic variety 
and $R$ is a birational isomorphism. 
The case of $\mathbb{CP}^1\times \mathbb{CP}^1$ has been recently discussed 
by Adler, Bobenko and Suris \cite{ABS1}, in relation with the classification of the 
so-called {\em quadrirational maps}. 
Remarkably, the resulting list of YB maps corresponds exactly to 
the five different ways of intersection of two conics. 
A classification for integrable equations on quad-graphs was also performed 
by the same authors in \cite{ABS2}, where the connection 
between the YB relation for maps and the {\em consistency} property 
for discrete equation on quad--graphs was addressed
(see concluding remarks of \cite{ABS2}). 
 
More recently, the connection between YB maps and integrable 
equations on quad--graphs has been investigated in \cite{VTS}. It was
shown that a systematic framework for studying such a connection is based
on the symmetry groups of the equations. The main observation is that the 
variables of certain YB maps can be chosen as invariants of the symmetry 
group of the corresponding lattice equation. Thus, it was demonstrated that 
YB maps can be obtained from  scalar integrable lattice equations
admitting a one-parameter group of symmetry transformations. 
Moreover, this idea was extended to scalar integrable lattice 
equations admitting a multi-parameter symmetry group by considering the 
extension of the equation on a multi-dimensional lattice,
as well as to multi-field integrable lattice equations.
 
In the present work, we continue the study on the connection between YB 
maps and  multi-field integrable lattice equations.
In section \ref{sec2} we present the explicit auto-transformations of 
integrable chains introduced by Adler and Yamilov in \cite{A-Y} and the
interpretation of these transformations as discrete equations on quad--graphs. 
Integrability aspects of such equations are discussed in section \ref{sec3}. 
There we stress the link between the three-dimensional
consistency property of lattice equations on quad-graphs and  
the braid type relation of the considered auto-transformations. 
Section \ref{sec4} deals with symmetry methods applied to 
lattice equations. The existence of a symmetry group of
transformations allows one to construct YB maps from
integrable lattice equations by considering  invariants of the 
transformation group as  YB variables. 
The procedure is applied to the list of lattice equations of 
section \ref{sec2} and the results are presented in section \ref{sec5}.
In section \ref{sec6} we derive YB maps from  $n$--vector extensions
 of the nonlinear Schr\"odinger system.
We conclude in section \ref{persp} by making comments 
on various perspectives of the subject. 

\section{Explicit auto-transformations of integrable chains as 2-field quad-graph equations}
\label{sec2}

Adler and Yamilov introduced in \cite{A-Y} explicit auto-transformations for integrable chains 
associated to systems of evolution partial differential equations of the form
\begin{equation}
u_t=u_{xx}+f(u,v,u_x)\,, \quad v_t=-v_{xx}+g(u,v,v_x) \,. \nonumber
\end{equation}
For example the nonlinear Schr\"odinger system 
\begin{equation}
u_t=u_{xx}+2\,u^2 \, v \,,\qquad v_t=-v_{xx}-2\, u\,v^2\,, \nonumber
\end{equation}
is related to the chain 
\begin{equation}
u_{jx}=-u_{j+1}-\alpha_{j}u_{j}- u_{j}^2 v_{j+1}\, , \quad 
v_{jx}= v_{j-1}+\alpha_{j-1}v_{j}+u_{j-1} v_{j}^2 \, . \nonumber
 \label{eq:nLSchain}
\end{equation}
For given $k\in\mathbb{Z}$, this chain admits 
(up to the exchange of two consecutive $\alpha$'s) 
an explicit auto--transformation 
$\mathbb{B}_k : (u_j,v_j) \, {\rightarrow}\, (\widetilde{u}_j, \widetilde{v}_j)$,
which preserves the values of the fields
$(u_j,v_j)$, $j\in \mathbb{Z}\setminus \lbrace k\rbrace $
and the parameters 
$\alpha_j$, $j\in \mathbb{Z}\setminus \lbrace k-1,k\rbrace$ 
and changes the remaining terms of the chain by
\begin{subequations} 
\begin{align}
& \widetilde{u}_k=u_k+\frac{\alpha_{k-1}-\alpha_k}{1-u_{k-1} v_{k+1}} \, u_{k-1}  \,,  \qquad
\widetilde{v}_k=v_k-\frac{\alpha_{k-1}-\alpha_k}{1-u_{k-1} v_{k+1}} \, v_{k+1} \, ,\\
&\widetilde{\alpha}_k=\alpha_{k-1}\, , \quad \, \tilde{\alpha}_{k-1}=\alpha_{k} \, .
\label{eq:AdlYamtilde}
\end{align}
\end{subequations} 
An interpretation of the above transformations can be given in terms of 
deformations of simple paths on the multidimensional orthogonal grid.
This interpretation follows naturally if one starts from the two dimensional grid situation
studied in \cite{PNC} and generalizes to higher dimension inspired by the multi-dimensional 
consistency \cite{N,ABS2} and the Cauchy problem on quad--graphs \cite{VseSa}.

\begin{figure}[h]
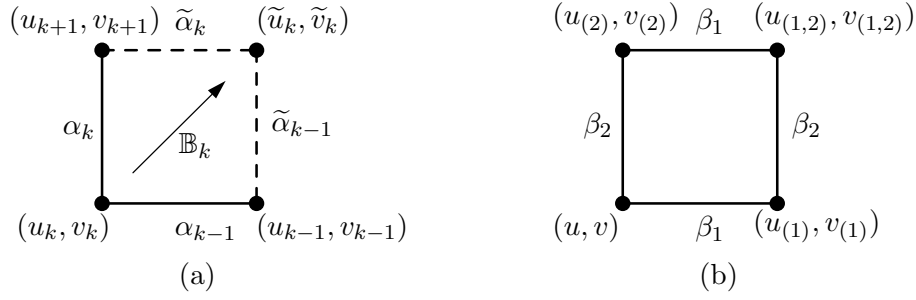

\begin{minipage}{15pc}
\centertexdraw{\setunitscale 0.4
\linewd 0.03 \arrowheadtype t:F \arrowheadsize l:0.25 w:0.15
\move (-1 -1)
\lvec (-1 1) \lpatt(0.15 0.15) \lvec (1 1) \lpatt(0.15 0.15) \lvec (1. -1) \lpatt() \lvec(-1 -1) 
\move(-0.6 -0.6) \linewd 0.02  \avec(0.6 0.6)  \linewd 0.03
\move(-1 -1) \fcir f:0.0 r:0.1 \move(-1 1) \fcir f:0.0 r:0.1 
\move(1 1) \fcir f:0.0 r:0.1 \move(1 -1) \fcir f:0.0 r:0.1    
\htext (-2.1 -1.5) {$(u_{k},v_{k})$} 
\htext (-2.2 1.2) {$(u_{k+1},v_{k+1})$}
\htext (1.0 1.2) {$(\widetilde{u}_k,\widetilde{v}_k)$}
\htext (1. -1.5) {$(u_{k-1},v_{k-1})$}
\htext (-0.05 1.2) {$\widetilde{\alpha}_{k}$}
\htext (-0.05 -1.5) {$\alpha_{k-1}$}
\htext (-1.5 -0.15) {$\alpha_{k}$}
\htext (1.2 -0.15) {$\widetilde{\alpha}_{k-1}$} 
\htext (0 -0.4) {$\mathbb{B}_k$} 
\htext (-3.5 1.8) {$\phantom{o}$} 
\htext (0.0 -2.15) {(a)}
}
\end{minipage}\hspace{6pc}%
\begin{minipage}{10pc}
\centertexdraw{\setunitscale 0.4
\linewd 0.03 \arrowheadtype t:F 
\move (-1 -1)
\lvec (-1 1) \lvec (1 1)  \lvec (1. -1) \lvec(-1 -1)
\move(-1 -1) \fcir f:0.0 r:0.1 \move(-1 1) \fcir f:0.0 r:0.1 
\move(1 1) \fcir f:0.0 r:0.1 \move(1 -1) \fcir f:0.0 r:0.1    
\htext (-1.9 -1.5) {$(u,v)$} 
\htext (-1.9 1.2) {$(u_{(2)},v_{(2)})$}
\htext (0.7 1.2) {$(u_{(1,2)},v_{(1,2)})$}
\htext (0.7 -1.5) {$(u_{(1)},v_{(1)})$}
\htext (-0.05 1.2) {$\beta_1$}
\htext (-0.05 -1.5) {$\beta_1$}
\htext (-1.5 -0.15) {$\beta_2$}
\htext (1.2 -0.15) {$\beta_2$}  
\htext (1.5 1.8) {$\phantom{o}$} 
\htext (0.0 -2.15) {(b)}
}
\end{minipage}
\caption{(a): An elementary flip representing the transformation $\mathbb{B}_k$. 
(b): The corresponding  notation of the fields and the parameters.}
\label{fig:Adlerflip}
\end{figure}

We consider the multidimensional orthogonal grid with set of vertices ${\mathbb Z}^m$,  
edges connecting neighboring vertices $(n_1,n_2,...,n_p,...,n_m)$, 
$(n_1,n_2,...,n_p+1,...,n_m)$ and a set of $m$ different complex parameters 
$\beta_1,\beta_2,\ldots,\beta_m$. For each $(n_1,n_2,...,n_p,...,n_m)\in {\mathbb Z}^m$
the edge between the vertices $(n_1,n_2,...,n_p,...,n_m)$ and $(n_1,n_2,...,n_p+1,...,n_m)$ 
is labeled by the value $\beta_p$, so parallel edges are labeled with the same  $\beta$-value.
A plane perpendicular to an edge at its midpoint is called a characteristic plane.  
We consider now a simple path on the grid defined by a sequence 
$\left(\kappa_j\right)$ of vertices with the property that 
each characteristic plane cuts one and only one edge of the path. 
This is an instance of a Cauchy path for quad--graph equations as was 
analyzed in \cite{VseSa}. At the vertex $\kappa_j$ we put the value of the 
2-field $(u_j,v_j)$ and we assume that the value of the parameter $\alpha_j$ 
of the chain coincides with the $\beta$--value of the edge $(\kappa_j$,$\kappa_{j+1})$.
Then, in the generic case when $\alpha_k \neq \alpha_{k-1}$, the transformation 
$\mathbb{B}_k$ can be thought as a local deformation replacing the  $(u_{k-1},v_{k-1})\stackrel{\alpha_{k-1}}{\rule[1pt]{1cm}{0.6pt}}(u_{k},v_{k})\stackrel{\alpha_{k}}
{\rule[1pt]{1cm}{0.6pt}}(u_{k+1},v_{k+1})$ angular segment of the path by 
$(u_{k-1},v_{k-1})\stackrel{\alpha_{k}}        {\rule[1pt]{1cm}{0.6pt}}(\tilde{u}_{k},\tilde{v}_{k})\stackrel{\alpha_{k-1}}{\rule[1pt]{1cm}{0.6pt}}(u_{k+1},v_{k+1})$ as in figure \ref{fig:Adlerflip}(a) i.e. flipping one angle of the square face to the opposite. 
Notice that if $\alpha_k=\alpha_{k-1}$, the transformation $\mathbb{B}_k$ 
reduces to the identity which is compatible with the fact that the corresponding 
two consecutive edges lie on the same line and there is no angle to flip.

Having this picture in mind we see that the basic building block is a square {\em plaquette} on a
2-dimensional plane (e.g. in the 1 and 2 directions) of the grid with equations (\ref{eq:AdlYamtilde})
rewritten in the form
\begin{equation}
u_{(1,2)}=u+\frac{\beta_{1}-\beta_2}{1-u_{(1)} v_{(2)}}\,u_{(1)} \,, \qquad
v_{(1,2)}=v-\frac{\beta_{1}-\beta_2}{1-u_{(1)} v_{(2)}}\,v_{(2)} \,, \label{eq:AdlYam8}
\end{equation}
where the subscript $i$ inside a parenthesis indicates the displacement of the fields $(u,v)$ 
in the $i$ direction by a unit step on the orthogonal grid (see figure \ref{fig:Adlerflip}(b) ). 
Following the above notation we list below the equations from \cite{A-Y} which are relevant 
in the present work. We denote them by ${\mathcal E}_a$,  $a=1,2,\ldots 5$, 
\begin{equation}
\mathcal{E}_{a}\,: \qquad 
u_{(1,2)}=E_{a}(u,u_{(1)},v_{(2)},\beta_{1},\beta_2) \, \, ,  \quad  
v_{(1,2)}=E_{a}(v,v_{(2)},u_{(1)},\beta_{2},\beta_1) \,,
\label{eq:AdlYamR}
\end{equation}
where the corresponding functions on the right-hand sides, evaluated at $(x,y,z,\beta_1,\beta_2)$, 
are given in table \ref{tab1} below.

\begin{table}[h]
\begin{tabular}{ll} \hline \\
$\displaystyle{E_{1}= x+(\beta_{1}-\beta_2)\,\frac{y}{1-y\, z} }$ &  
$\quad\,\,\,\,\, \displaystyle{E_{4}=  x + (\frac{1}{\beta_1}-\frac{1}{\beta_2})\, 
\frac{\beta_1 \, \beta_2\, (x + z)\, (y-x) - \epsilon}{\beta_2\,(x+z) + \beta_1\,(y-x)} } $ \hfill $\qquad\qquad\qquad$
\\ \\
$\displaystyle{E_{2}= x + (\beta_2-\beta_1)\, \frac{y-x}{\beta_2 + y\, z} }$ & \\ \\
$\displaystyle{E_{3}= x + (\beta_1-\beta_2)\, \frac{x + z}{y+z -\beta_1} }$ &  
$\quad\,\,\,\,\, \displaystyle{E_{5}=  \frac{(\beta_1-\beta_2)\, y \, z + 
x \, (\gamma_2 \, y+ \gamma_1 \, z)}{ (\beta_1- \beta_2)\,  x + \gamma_1 \, y + \gamma_2 \, z} }$ \\ \\ \hline 
\end{tabular}
\caption{ \label{tab1} The list of functions $E_{a}(x,y,z,\beta_1,\beta_2)$. 
The additional parameters $\gamma_i$, appearing in function $E_5$, depend on
the parameters $\beta_i$ through the relations ${\gamma_i}^2 - {\beta_i}^2 = \delta$, 
where $\delta$ can be fixed to be $0$ or $1$. In function $E_4$, $\epsilon$ is a free parameter.} 
\end{table}
\begin{rem}
System $ {\mathcal E}_1$ corresponds to equations (\ref{eq:AdlYam8}). System ${\mathcal E}_{2}$ 
is related to the chain of the derivative nonlinear Schr\"odinger system 
\cite{ChenLeeLiu} deformed by its first symmetry.
System ${\mathcal E}_{3}$ is related to the chain 
associated to the Kaup system \cite{Kaup}.
The remaining lattice equations $ {\mathcal E}_{4} $ and $ {\mathcal E}_{5}$ 
are related to chains corresponding to  the Landau-Lifshits system
\begin{equation}u_t=u_{xx}- \frac{2}{u+v}\big(u_x^2+P(u)\big)+\frac{1}{2} P'(u)\, , \,  \,\quad 
v_t=-v_{xx}+ \frac{2}{u+v}\big(v_x^2+P(-v)\big)+\frac{1}{2} P'(-v)\,, \nonumber
\end{equation}
for the specific cases where $P(u) \equiv \epsilon$ (system ${\mathcal E}_{4}$) 
and $P(u) \equiv \delta u^2$ (system ${\mathcal E}_{5}$). 
\end{rem}

\section{Braid relation and the three-dimensional consistency property}
\label{sec3}

Key properties of the transformations $\mathbb{B}_j$ are the relations $\mathbb{B}_j^2=1\, , \, (\mathbb{B}_j \, \mathbb{B}_{j+1})^3=1 $  and $ \mathbb{B}_j \, \mathbb{B}_i=\mathbb{B}_i \mathbb{B}_j $, for $|i-j|>1$.
The first one means that each transformation $\mathbb{B}_j$ is an involution. 
The second one, in view of the first, yields the following braid-type relation 
\begin{equation}
\mathbb{B}_{j+1}\,\mathbb{B}_j\,\mathbb{B}_{j+1}=\mathbb{B}_j\,\mathbb{B}_{j+1}\,\mathbb{B}_j.
\label{eq:Braid}
\end{equation}

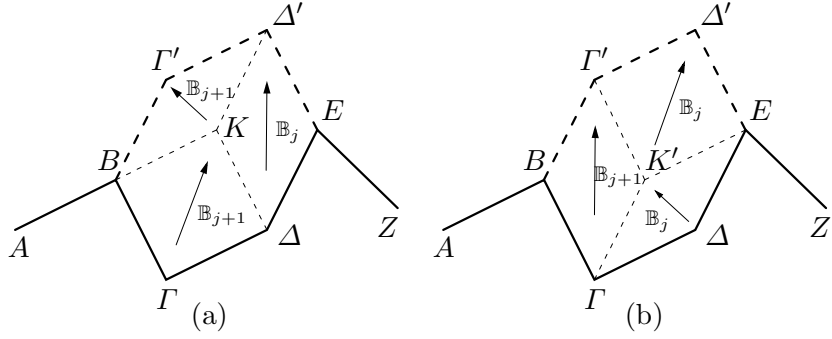
\begin{figure}[h] 
\setlength{\unitlength}{0.00029in}
\begingroup\makeatletter\ifx\SetFigFont\undefined%
\gdef\SetFigFont#1#2#3#4#5{%
  \reset@font\fontsize{#1}{#2pt}%
  \fontfamily{#3}\fontseries{#4}\fontshape{#5}%
  \selectfont}%
\fi\endgroup%
{\renewcommand{\dashlinestretch}{30}
\begin{picture}(14360,6000)(-3500,-900)
\thicklines
\path(7683,933)(9483,1833)(10383,33)
(12183,933)(13083,2733)(14519,1326)
\dashline{200.000}(9483,1833)(10383,3633)(12183,4533)
\dashline{200.000}(12183,4533)(13083,2733)
\thinlines
\dashline{100.000}(10383,33)(11283,1833)(13083,2733)
\dashline{100.000}(11283,1833)(10383,3633)
\thicklines
\path(33,933)(1833,1833)(2733,33)
(4533,933)(5433,2733)(6869,1326)
\thinlines
\dashline{100.000}(1833,1833)(3633,2733)(4533,933)
\thicklines
\dashline{200.000}(1833,1833)(2733,3633)(4533,4533)
\thinlines
\dashline{100.000}(3633,2733)(4533,4533)
\thicklines
\dashline{200.000}(4533,4533)(5433,2733)
\thinlines
\path(2913,663)(3453,2103)
\blacken\thicklines
\path(3440.711,1963.428)(3453.000,2103.000)(3370.486,1989.763)(3440.711,1963.428)
\thinlines
\path(4519,2005)(4533,3543)
\blacken\thicklines
\path(4569.270,3407.664)(4533.000,3543.000)(4494.273,3408.347)(4569.270,3407.664)
\thinlines
\path(10383,1203)(10397,2741)
\blacken\thicklines
\path(10433.270,2605.664)(10397.000,2741.000)(10358.273,2606.347)(10433.270,2605.664)
\thinlines
\path(12067,1084)(11472,1630)
\blacken\thinlines
\path(11596.822,1566.354)(11472.000,1630.000)(11546.113,1511.094)(11596.822,1566.354)
\path(11463,2463)(11976,3912)
\blacken\thicklines
\path(11966.295,3772.225)(11976.000,3912.000)(11895.595,3797.255)(11966.295,3772.225)
\thinlines
\path(3453,2913)(2858,3459)
\blacken\thicklines
\path(2982.822,3395.354)(2858.000,3459.000)(2932.113,3340.094)(2982.822,3395.354)
\put(7500,450){\makebox(0,0)[lb]{{\SetFigFont{11}{8}{\rmdefault}{\mddefault}{\updefault}$A$}}}
\put(9100,2013){\makebox(0,0)[lb]{{\SetFigFont{11}{10}{\rmdefault}{\mddefault}{\updefault}$B$}}}
\put(10203,-500){\makebox(0,0)[lb]{{\SetFigFont{11}{11}{\rmdefault}{\mddefault}{\updefault}$\varGamma$}}}
\put(12363,663){\makebox(0,0)[lb]{{\SetFigFont{11}{11}{\rmdefault}{\mddefault}{\updefault}$\varDelta$}}}
\put(13150,2913){\makebox(0,0)[lb]{{\SetFigFont{11}{11}{\rmdefault}{\mddefault}{\updefault}$E$}}}
\put(14220,800){\makebox(0,0)[lb]{{\SetFigFont{11}{11}{\rmdefault}{\mddefault}{\updefault}$Z$}}}
\put(-100,450){\makebox(0,0)[lb]{{\SetFigFont{11}{11}{\rmdefault}{\mddefault}{\updefault}$A$}}}
\put(1500,2013){\makebox(0,0)[lb]{{\SetFigFont{11}{11}{\rmdefault}{\mddefault}{\updefault}$B$}}}
\put(2553,-500){\makebox(0,0)[lb]{{\SetFigFont{11}{12}{\rmdefault}{\mddefault}{\updefault}$\varGamma$}}}
\put(4713,663){\makebox(0,0)[lb]{{\SetFigFont{11}{11}{\rmdefault}{\mddefault}{\updefault}$\varDelta$}}}
\put(5503,2913){\makebox(0,0)[lb]{{\SetFigFont{11}{11}{\rmdefault}{\mddefault}{\updefault}$E$}}}
\put(6500,800){\makebox(0,0)[lb]{{\SetFigFont{11}{11}{\rmdefault}{\mddefault}{\updefault}$Z$}}}
\put(4623,4623){\makebox(0,0)[lb]{{\SetFigFont{11}{11}{\rmdefault}{\mddefault}{\updefault}${\varDelta}'$}}}
\put(2463,3813){\makebox(0,0)[lb]{{\SetFigFont{11}{14.4}{\rmdefault}{\mddefault}{\updefault}${\varGamma}'$}}}
\put(3743,2600){\makebox(0,0)[lb]{{\SetFigFont{11}{14.4}{\rmdefault}{\mddefault}{\updefault}$K$}}}
\put(11283,2013){\makebox(0,0)[lb]{{\SetFigFont{11}{14.4}{\rmdefault}{\mddefault}{\updefault}$K'$}}}
\put(12273,4623){\makebox(0,0)[lb]{{\SetFigFont{11}{14.4}{\rmdefault}{\mddefault}{\updefault}${\varDelta}'$}}}
\put(10113,3813){\makebox(0,0)[lb]{{\SetFigFont{11}{14.4}{\rmdefault}{\mddefault}{\updefault}${\varGamma}'$}}}
\put(3100,3250){\makebox(0,0)[lb]{{\SetFigFont{8}{14.4}{\rmdefault}{\mddefault}{\updefault}${\mathbb{B}_{j+1}}$}}}
\put(4700,2500){\makebox(0,0)[lb]{{\SetFigFont{8}{14.4}{\rmdefault}{\mddefault}{\updefault}${\mathbb{B}_{j}}$}}}
\put(3300,1000){\makebox(0,0)[lb]{{\SetFigFont{8}{14.4}{\rmdefault}{\mddefault}{\updefault}${\mathbb{B}_{j+1}}$}}}
\put(11300,900){\makebox(0,0)[lb]{{\SetFigFont{8}{14.4}{\rmdefault}{\mddefault}{\updefault}${\mathbb{B}_{j}}$}}}
\put(11900,2900){\makebox(0,0)[lb]{{\SetFigFont{8}{14.4}{\rmdefault}{\mddefault}{\updefault}${\mathbb{B}_{j}}$}}}
\put(10380,1700){\makebox(0,0)[lb]{{\SetFigFont{8}{14.4}{\rmdefault}{\mddefault}{\updefault}${\mathbb{B}_{j+1}}$}}}
\put(3000,-900){\makebox(0,0)[lb]{{\SetFigFont{11}{14.4}{\rmdefault}{\mddefault}{\updefault} (a) }}}
\put(10750,-900){\makebox(0,0)[lb]{{\SetFigFont{11}{14.4}{\rmdefault}{\mddefault}{\updefault} (b) }}}
\end{picture}
}
\caption{The braid relation as the consistency property for the deformations of a path on the multidimensional grid.}
\label{fig:flip}
\end{figure}

Let us now explain the implication of the braid relation (\ref{eq:Braid}) in terms 
of deformations of paths on the grid. 
By applying the composite transformations appearing on each side of equation (\ref{eq:Braid}) to a 
path $AB \varGamma \varDelta EZ$, as shown in figure \ref{fig:flip}, we get two
deformations. 
Consider first the left hand side of the braid relation: $\mathbb{B}_{j+1}$ flips the $\varGamma$ corner to $K$, then 
$\mathbb{B}_{j} $ flips the $\varDelta$ corner to $\varDelta'$ and finally  $\mathbb{B}_{j+1}$ 
flips the $K$ corner to $\varGamma'$.
For the right hand side, $\mathbb{B}_{j}$ flips the $\varDelta$ corner to $K'$ then 
$\mathbb{B}_{j+1} $ flips the $\varGamma$ corner to $\varGamma'$  and finally 
$\mathbb{B}_{j}$ flips the $K'$ corner to $\varDelta'$.
Thus the original path $ AB \varGamma \varDelta EZ$ is transformed to the same path 
$AB\varGamma'\varDelta'EZ$ through two different deformations according to figure \ref{fig:flip}, viz.  
\begin{align}
{\rm{(a)}}\,:&\qquad AB \varGamma \varDelta EZ \stackrel{\mathbb{B}_{j+1}}{\longrightarrow} ABK \varDelta EZ \stackrel{\mathbb{B}_{j}}{\longrightarrow}ABK\varDelta'EZ  \stackrel{\mathbb{B}_{j+1}}{\longrightarrow}AB\varGamma'\varDelta'EZ \,, \nonumber \\
{\rm{(b)}}\,:&\qquad AB \varGamma \varDelta EZ\stackrel{\mathbb{B}_{j}}{\longrightarrow}
AB \varGamma K'EZ\stackrel{\mathbb{B}_{j+1}}{\longrightarrow}   
AB\varGamma'K'EZ\stackrel{\mathbb{B}{j}}{\longrightarrow}AB\varGamma'\varDelta'EZ \,. \nonumber
\end{align}

The braid relation (\ref{eq:Braid}) guarantees that the values at the vertices 
$\varGamma'$, $\varDelta'$ are uniquely determined in terms of the values at the vertices 
$B$, $\varGamma$, $\varDelta$, $E$ i.e. independent of the way we deform the path in order to 
arrive at these vertices. 
This is an instance of the three-dimensional consistency for quad-graph equations 
(see \cite{N, BS, ABS2}, as well as \cite{VTS}). The consistency property can be checked for an initial 
value configuration as in figure \ref{fig:flip}, with values assigned on the vertices 
$B$, $\varGamma$, $\varDelta$, $E$.
This configuration is adapted to the cubic interpretation of the YB property \cite{Suris}, 
cf also \cite{VTS}.

\section{Invariants of symmetry groups and Yang--Baxter variables}
\label{sec4}

We start with the building block of equations on quad-graphs consisting of a system of 
algebraic relations of the form
\begin{equation}
\mathcal{B}^i(f,f_{(1)},f_{(2)},f_{(1,2)};\beta_1,\beta_2)=0\,, \quad i=1,\ldots \ell \,,\label{eq:quad} 
\end{equation}
where $f=(f^1,f^2,\ldots,f^{\ell})$ and $f^i:\mathbb{Z}^2 \rightarrow \mathbb{C}$ (or $\mathbb{CP}^1$), 
$i=1,2,\ldots,\ell$. We use upper indices to denote the components of the fields, since we have reserved lower
indices to indicate their shifted values.
Consider a one-parameter group of transformations $G$ acting on   
the domain of the dependent variables of a lattice equation (\ref{eq:quad}), i.e.
\begin{equation}
G: f^i \mapsto \Phi^i(n_1,n_2,f;\varepsilon) \,, \quad \varepsilon \in \mathbb{C}\,, \qquad i=1,\ldots \ell\,, \nonumber
\end{equation}
where $n_1,n_2 \in \mathbb{Z}$ denote the site variables of the lattice.
Let ${\rm J}^{(k)}$ denote the lattice jet space with coordinates $(f,f_{J})$, where by
$f_{J}$ we mean the forward shifted values of $f$, indexed by all unordered (symmetric) multi-indices 
$J=(j^1,j^2,\ldots j^k)$, $1\leq j^r \leq 2$, of order $k=\# J$. 
The prolongation of the group action of $G$ on ${\rm J}^{(k)}$ is
\begin{equation}
G^{(k)}: (f^i, f^i_J) \mapsto (\Phi^i(n_1,n_2,f;\varepsilon), 
\Phi^i_{J}(n_1,n_2,f;\varepsilon))\,,\quad i=1,\ldots \ell \,,\label{eq:Gpr} 
\end{equation}
where $\Phi^i_{(1)}(n_1,n_2,f;\varepsilon)=\Phi^i(n_1+1,n_2,f_{(1)};\varepsilon)$, etc.

The infinitesimal generator of the group action of $G$ on the space of the dependent variables is given by the vector field
\begin{equation}
\mathbf{v} = \sum_{i=1}^n Q^i(n_1,n_2,f)\,\partial_{f^i}\,, \quad \mbox{where}\quad  
Q^i(n_1,n_2,f)=\left. 
\frac{\rm d \phantom{\varepsilon}}{{\rm d}\, \varepsilon}\,\Phi^i(n_1,n_2,f;\varepsilon) \right|_{\varepsilon=0} \,.
\nonumber
\end{equation}
There is a one-to-one correspondence between connected groups of transformations and their associated 
infinitesimal generators, since the group action is reconstructed by the flow of the vector field $\mathbf{v}$ 
by exponentiation
\begin{equation}
\Phi^i(n_1,n_2,f;\varepsilon)=\exp (\varepsilon \, \mathbf{v}) f^i\,, \quad i=1,\ldots n \,.  \nonumber
\end{equation}
The infinitesimal generator of the action of $G^{(k)}$ on ${\rm J}^{(k)}$ is the associated $k^{\rm th}$ order
prolonged vector field
\begin{equation}
\mathbf{v}^{(k)} = \sum_{i=1}^n \sum_{\# J=j=0}^k Q_J^i(n_1,n_2,f)\,\partial_{f^i_J}\,. \nonumber
\end{equation}
The transformation $G$ is a Lie-point symmetry of the lattice equations (\ref{eq:quad}), if it transforms any solution of 
(\ref{eq:quad}) to another solution of the same equations. Equivalently, $G$ is a symmetry of equations 
(\ref{eq:quad}), if the equations are not affected by the transformation (\ref{eq:Gpr}).
The infinitesimal criterion for a connected group of transformations $G$ to be a symmetry of equations 
(\ref{eq:quad}) is 
\begin{equation}
\mathbf{{v}}^{(2)} \big(\mathcal{B}^i(f,f_{(1)},f_{(2)},f_{(1,2)};\beta_1,\beta_2)\big) =0\,, \label{eq:infinv}
\end{equation}
$i=1,2\ldots,n$. Equations (\ref{eq:infinv}) should hold for all solutions of equations (\ref{eq:quad}), 
and thus the latter equations and their consequences should be taken into account. 
Equations (\ref{eq:infinv}) determine the most general infinitesimal Lie point
symmetry of the system (\ref{eq:quad}). The resulting set of infinitesimal
generators forms a Lie algebra $\mathfrak{g}$ from which the corresponding 
Lie point symmetry group $G$ can be found by exponentiating the given vector fields.

A {\em lattice invariant} under the action of $G$ is a function $I:{\rm J}^{(k)}\rightarrow \mathbb{C}$ 
which satisfies $I(g^{(k)} \cdot (f, f_J)) = I( f,f_J)$ for all $g \in G$ and all 
$( f,f_j) \in {\rm J}^{(k)}$.
For connected groups of transformations, a necessary and sufficient condition for a function 
$I$ to be invariant under the action of $G$, is the annihilation 
of $I$ by all prolonged infinitesimal generators, i.e.
\begin{equation}
\mathbf{{v}}^{(k)}(I) = 0 \,. \label{eq:pde}
\end{equation}
for all $\mathbf{{v}} \in \mathfrak{g}$.

\begin{figure}[h]
\setlength{\unitlength}{0.0005in}
\begingroup\makeatletter\ifx\SetFigFont\undefined%
\gdef\SetFigFont#1#2#3#4#5{%
  \reset@font\fontsize{#1}{#2pt}%
  \fontfamily{#3}\fontseries{#4}\fontshape{#5}%
  \selectfont}%
\fi\endgroup%
{\renewcommand{\dashlinestretch}{30}
\begin{picture}(7224,4000)(-2750,-700)
\thinlines
\path(912,2712)(912,912)(2712,912)
	(2712,2712)(912,2712)
\path(912,2712)(12,1812)(12,12)(912,912)
\path(12,12)(1812,12)(2712,912)
\path(4512,1812)(4512,12)(6312,12)
	(6312,1812)(4512,1812)
\path(6312,12)(7212,912)(7212,2712)(6312,1812)
\path(4512,1812)(5412,2712)(7212,2712)
\path(7212,1812)(6312,912)
\blacken\path(6375.640,1018.066)(6312.000,912.000)(6418.066,975.640)(6375.640,1018.066)
\path(6312,912)(4512,912)
\blacken\path(4632.000,942.000)(4512.000,912.000)(4632.000,882.000)(4632.000,942.000)
\path(2712,1812)(912,1812)
\blacken\path(1032.000,1842.000)(912.000,1812.000)(1032.000,1782.000)(1032.000,1842.000)
\path(6762,2262)(4962,2262)
\blacken\path(5082.000,2292.000)(4962.000,2262.000)(5082.000,2232.000)(5082.000,2292.000)
\path(2262,462)(462,462)
\blacken\path(582.000,492.000)(462.000,462.000)(582.000,432.000)(582.000,492.000)
\path(912,1812)(12,912)
\blacken\path(75.640,1018.066)(12.000,912.000)(118.066,975.640)(75.640,1018.066)
\path(1812,1880)(1812,2712)
\blacken\path(1842.000,2592.000)(1812.000,2712.000)(1782.000,2592.000)(1842.000,2592.000)
\path(1812,912)(1812,1744)
\path(462,1452)(462,2262)
\blacken\path(492.000,2142.000)(462.000,2262.000)(432.000,2142.000)(492.000,2142.000)
\path(462,462)(462,1294)
\path(912,12)(1317,417)
\path(1407,507)(1812,912)
\blacken\path(1748.360,805.934)(1812.000,912.000)(1705.934,848.360)(1748.360,805.934)
\path(6762,462)(6762,1294)
\path(6762,1452)(6762,2262)
\blacken\path(6792.000,2142.000)(6762.000,2262.000)(6732.000,2142.000)(6792.000,2142.000)
\path(5412,980)(5412,1812)
\blacken\path(5442.000,1692.000)(5412.000,1812.000)(5382.000,1692.000)(5442.000,1692.000)
\path(5412,12)(5412,844)
\path(5412,1812)(5817,2217)
\path(5907,2307)(6312,2712)
\blacken\path(6248.360,2605.934)(6312.000,2712.000)(6205.934,2648.360)(6248.360,2605.934)
\put(6550,950){\makebox(0,0)[lb]{{\SetFigFont{10}{14.4}{\rmdefault}{\mddefault}{\updefault}$R^{(2,3)}$}}}
\put(5450,600){\makebox(0,0)[lb]{{\SetFigFont{10}{14.4}{\rmdefault}{\mddefault}{\updefault}$R^{(1,3)}$}}}
\put(5800,1950){\makebox(0,0)[lb]{{\SetFigFont{10}{14.4}{\rmdefault}{\mddefault}{\updefault}$R^{(1,2)}$}}}
\put(1900,1500){\makebox(0,0)[lb]{{\SetFigFont{10}{14.4}{\rmdefault}{\mddefault}{\updefault}$R^{(1,3)}$}}}
\put(1300,170){\makebox(0,0)[lb]{{\SetFigFont{10}{14.4}{\rmdefault}{\mddefault}{\updefault}$R^{(1,2)}$}}}
\put(230,950){\makebox(0,0)[lb]{{\SetFigFont{10}{14.4}{\rmdefault}{\mddefault}{\updefault}$R^{\,(2,3)}$}}}
\put(800,-300){\makebox(0,0)[lb]{{\SetFigFont{10}{14.4}{\rmdefault}{\mddefault}{\updefault}$x$}}}
\put(1900,1000){\makebox(0,0)[lb]{{\SetFigFont{10}{14.4}{\rmdefault}{\mddefault}{\updefault}$x'$}}}
\put(1900,2800){\makebox(0,0)[lb]{{\SetFigFont{10}{14.4}{\rmdefault}{\mddefault}{\updefault}$x''$}}}
\put(2400,300){\makebox(0,0)[lb]{{\SetFigFont{10}{14.4}{\rmdefault}{\mddefault}{\updefault}$y$}}}
\put(250,400){\makebox(0,0)[lb]{{\SetFigFont{10}{14.4}{\rmdefault}{\mddefault}{\updefault}$y'$}}}
\put(2800,1750){\makebox(0,0)[lb]{{\SetFigFont{10}{14.4}{\rmdefault}{\mddefault}{\updefault}$z$}}}
\put(700,1750){\makebox(0,0)[lb]{{\SetFigFont{10}{14.4}{\rmdefault}{\mddefault}{\updefault}$z'$}}}
\put(400,2400){\makebox(0,0)[lb]{{\SetFigFont{10}{14.4}{\rmdefault}{\mddefault}{\updefault}$y''$}}}
\put(-300,800){\makebox(0,0)[lb]{{\SetFigFont{10}{14.4}{\rmdefault}{\mddefault}{\updefault}$z''$}}}
\put(5300,-300){\makebox(0,0)[lb]{{\SetFigFont{10}{14.4}{\rmdefault}{\mddefault}{\updefault}$x$}}}
\put(6900,300){\makebox(0,0)[lb]{{\SetFigFont{10}{14.4}{\rmdefault}{\mddefault}{\updefault}$y$}}}
\put(7300,1750){\makebox(0,0)[lb]{{\SetFigFont{10}{14.4}{\rmdefault}{\mddefault}{\updefault}$z$}}}
\put(6700,2300){\makebox(0,0)[lb]{{\SetFigFont{10}{14.4}{\rmdefault}{\mddefault}{\updefault}$y^\ast$}}}
\put(4800,2300){\makebox(0,0)[lb]{{\SetFigFont{10}{14.4}{\rmdefault}{\mddefault}{\updefault}$y^{\ast\ast}$}}}
\put(5300,1900){\makebox(0,0)[lb]{{\SetFigFont{10}{14.4}{\rmdefault}{\mddefault}{\updefault}$x^\ast$}}}
\put(6300,2800){\makebox(0,0)[lb]{{\SetFigFont{10}{14.4}{\rmdefault}{\mddefault}{\updefault}$x^{\ast\ast}$}}}
\put(6050,1000){\makebox(0,0)[lb]{{\SetFigFont{10}{14.4}{\rmdefault}{\mddefault}{\updefault}$z^\ast$}}}
\put(4150,900){\makebox(0,0)[lb]{{\SetFigFont{10}{14.4}{\rmdefault}{\mddefault}{\updefault}$z^{\ast\ast}$}}}
\put(1300,-700){\makebox(0,0)[lb]{{\SetFigFont{10}{14.4}{\rmdefault}{\mddefault}{\updefault}(a)}}}
\put(5700,-700){\makebox(0,0)[lb]{{\SetFigFont{10}{14.4}{\rmdefault}{\mddefault}{\updefault}(b)}}}
\end{picture}
}
\caption{Three dimensional representation of the YB relation}
\label{fig:3DYB}
\end{figure}
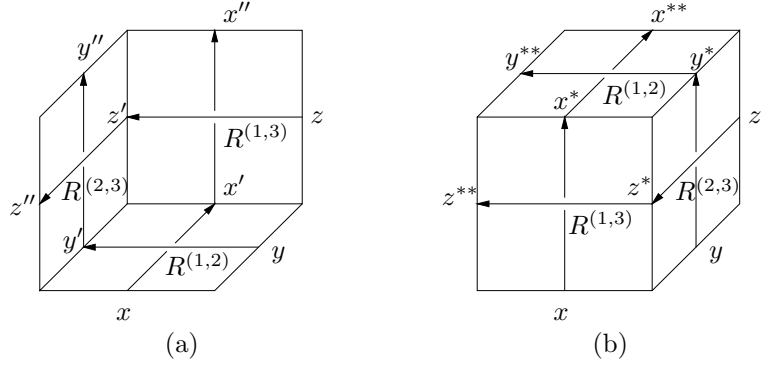

Consider a system of lattice equations of the form (\ref{eq:quad}) ($\ell=2$) and suppose that it admits 
a two-parameter group of symmetry transformations $G$ generated by two 
vector fields $\{\mathbf{v}_1,\mathbf{v}_2 \}$.
The invariants under the action of $G$ on ${\rm J}^{(2)}$ are obtained from the general solution 
of the system of first order partial differential equations (\ref{eq:pde}), by using the method of characteristics.
Using these invariant functions $I^i$, we assign now to the edges of an elementary quadrilateral 
the following YB variables
\begin{align}
x^i=I^i(f,f_{(1)}), \,\,\, y^i=I^i(f_{(2)},f),\,\,\, p^i=I^i(f_{(2)},f_{(1,2)}),\,\,\, q^i=I^i(f_{(1,2)},f_{(1)})\,.
\label{eq:edgeinv}
\end{align}
The invariants $x$, $y$, $p$, $q$, where $x=(x^1,x^2)$ etc, are not independent; 
there exist two relations among them
\begin{equation}
\mathcal{F}^i(x,y,u,v)=0 \,,\qquad i=1,2\,, \label{eq:rel}
\end{equation}
following from the fact that the space of $G$-orbits on ${\rm J}^{(2)}$ is six-dimensional.
Moreover, since $G$ is a symmetry group of the system of equations (\ref{eq:quad}), it can be written
in terms of the variables (\ref{eq:edgeinv}), yielding two more relations among the invariants $x$, $y$, $p$, $q$, i.e.
\begin{equation}
\mathcal{D}^i(x,y,p,q; \beta_1,\beta_2)=0 \,,\qquad i=1,2 \,. \label{eq:invlattice}
\end{equation}
Solving the system of equations (\ref{eq:invlattice}), (\ref{eq:rel}) for $p,q$
in terms of $x,y$ and assuming that the solution is unique, we obtain a two-parameter
family of maps $R(\beta_1,\beta_2) \, :\, (x,y)\mapsto (p,q)$.
The resulting map $R$ satisfies automatically the parameter dependent YB relation
\begin{equation}
R^{(2,3)} (\beta_2,\beta_3) \, R^{(1,3)} (\beta_1,\beta_3) \, R^{(1,2)}
(\beta_1,\beta_2) = R^{(1,2)} (\beta_1,\beta_2) \, R^{(1,3)}
(\beta_1,\beta_3) \, R^{(2,3)} (\beta_2,\beta_3) \,.
\label{eq:YBcom}
\end{equation}
This fact follows immediately from the interpretation of the braid relation (\ref{eq:Braid}) and
the initial value problem discussed at the end of the previous section, in terms
of YB variables (see figure \ref{fig:3DYB}). More precisely, the LHS and RHS of the braid relation 
(\ref{eq:Braid}) correspond to
the composition of mappings (a), (b) respectively, as follows
\begin{align}
{\rm (a)}:& \quad (x,y,z) & 
\stackrel{R^{(1,2)}}{\longrightarrow} &&({x'},{y'},z) \,&& 
\stackrel{R^{(1,3)}}{\longrightarrow} &&({x''},{y'},z')\,\, && 
\stackrel{R^{(2,3)}}{\longrightarrow} &&({x''},{y''},{z''}) \,\,\,\, , \nonumber
\\
{\rm (b)}:& \quad (x,y,z) & 
\stackrel{R^{(2,3)}}{\longrightarrow} &&(x,y^{\ast},z^{\ast}) && 
\stackrel{R^{(1,3)}}{\longrightarrow} &&(x^{\ast},y^{\ast},z^{\ast\ast}) && 
\stackrel{R^{(1,2)}}{\longrightarrow} &&(x^{\ast\ast},y^{\ast\ast},z^{\ast\ast}) \,. \nonumber
\end{align}

The preceding considerations can be casted in the form of the following:
\begin{prop} \label{propos}
If the system of discrete equations $\mathcal{B}^i=0$, $i=1,2$, satisfy the three dimensional 
consistency property and admits a two--parameter group of symmetry transformations, 
then the map $R(\beta_1,\beta_2)$ which relates the lattice invariants (\ref{eq:edgeinv}), 
satisfies the parameter dependent YB relation.
\end{prop}

\begin{figure}[h] 
\setlength{\unitlength}{0.0004in}
\begingroup\makeatletter\ifx\SetFigFont\undefined%
\gdef\SetFigFont#1#2#3#4#5{%
  \reset@font\fontsize{#1}{#2pt}%
  \fontfamily{#3}\fontseries{#4}\fontshape{#5}%
  \selectfont}%
\fi\endgroup%
{\renewcommand{\dashlinestretch}{30}
\begin{picture}(3246,4500)(-6500,-10)
\thinlines
\thinlines
\thinlines
\path(180,1530)(1530,180)
\path(1530,1080)(1530,180)
\path(1530,180)(2880,1530)
\texture{44000000 aaaaaa aa000000 8a888a 88000000 aaaaaa aa000000 888888 
	88000000 aaaaaa aa000000 8a8a8a 8a000000 aaaaaa aa000000 888888 
	88000000 aaaaaa aa000000 8a888a 88000000 aaaaaa aa000000 888888 
	88000000 aaaaaa aa000000 8a8a8a 8a000000 aaaaaa aa000000 888888 }
\shade\path(1530,1080)(2880,2430)(2880,1530)(1530,1080)
\path(1530,1080)(2880,2430)(2880,1530)(1530,1080) \thinlines
\path(1530,2880)(2430,1980)
\shade\path(1530,3780)(2880,2430)(1530,2880)(1530,3780)
\path(1530,3780)(2880,2430)(1530,2880)(1530,3780) \thinlines
\path(1530,2880)(630,1980)
\shade\path(180,2430)(1530,2880)(1530,3780)(180,2430) \thinlines
\path(180,2430)(1530,2880)(1530,3780)(180,2430)
\shade\path(180,2385)(180,1530)(1530,1080)(180,2430) \thinlines
\path(180,2385)(180,1530)(1530,1080)(180,2430)
\blacken\path(2176.066,2276.360)(2070.000,2340.000)(2133.640,2233.934)(2129.397,2280.603)(2176.066,2276.360)
\thinlines
\blacken\thinlines
\path(916.066,1736.360)(810.000,1800.000)(873.640,1693.934)(869.397,1740.603)(916.066,1736.360)
\thinlines
\blacken\thinlines
\path(843.329,1340.513)(720.000,1350.000)(824.355,1283.592)(799.689,1323.437)(843.329,1340.513)
\thinlines
\path(720,1350)(855,1305)
\blacken\thinlines
\path(926.360,2233.934)(990.000,2340.000)(883.934,2276.360)(930.603,2280.603)(926.360,2233.934)
\thinlines
\blacken\thinlines
\path(836.360,3043.934)(900.000,3150.000)(793.934,3086.360)(840.603,3090.603)(836.360,3043.934)
\thinlines
\blacken\thinlines
\path(2193.329,2690.513)(2070.000,2700.000)(2174.355,2633.592)(2149.689,2673.437)(2193.329,2690.513)
\thinlines
\path(2070,2700)(2205,2655)
\blacken\thinlines
\path(2186.360,1693.934)(2250.000,1800.000)(2143.934,1736.360)(2190.603,1740.603)(2186.360,1693.934)
\thinlines
\blacken\thinlines
\path(2276.360,883.934)(2340.000,990.000)(2233.934,926.360)(2280.603,930.603)(2276.360,883.934)
\thinlines
\dashline{100.000}(630,1980)(180,1530)
\dashline{100.000}(2430,1980)(2880,1530)
\blacken\path(885.645,2633.592)(990.000,2700.000)(866.671,2690.513)(910.311,2673.437)(885.645,2633.592)
\thinlines
\path(990,2700)(855,2655)
\blacken\thinlines
\path(2266.066,3086.360)(2160.000,3150.000)(2223.640,3043.934)(2219.397,3090.603)(2266.066,3086.360)
\thinlines
\blacken\thinlines
\path(826.066,926.360)(720.000,990.000)(783.640,883.934)(779.397,930.603)(826.066,926.360)
\thinlines
\blacken\thinlines
\path(2235.645,1283.592)(2340.000,1350.000)(2216.671,1340.513)(2260.311,1323.437)(2235.645,1283.592)
\thinlines
\path(2340,1350)(2205,1305)
\put(180,1530){\blacken\ellipse{90}{90}}
\put(180,1530){\ellipse{90}{90}}
\put(1530,180){\blacken\ellipse{90}{90}}
\put(1530,180){\ellipse{90}{90}}
\put(2880,1530){\blacken\ellipse{90}{90}}
\put(2880,1530){\ellipse{90}{90}}
\put(1530,1080){\blacken\ellipse{90}{90}}
\put(1530,1080){\ellipse{90}{90}}
\put(1530,2880){\blacken\ellipse{90}{90}}
\put(1530,2880){\ellipse{90}{90}}
\put(180,2430){\blacken\ellipse{90}{90}}
\put(180,2430){\ellipse{90}{90}}
\put(1530,3780){\blacken\ellipse{90}{90}}
\put(1530,3780){\ellipse{90}{90}}
\put(2880,2430){\blacken\ellipse{90}{90}}
\put(2880,2430){\ellipse{90}{90}}
\put(1300,2300){\makebox(0,0)[lb]{{\SetFigFont{11}{14.4}{\rmdefault}{\mddefault}{\updefault}$v_{(2)}$}}}
\put(1400,-200){\makebox(0,0)[lb]{{\SetFigFont{11}{14.4}{\rmdefault}{\mddefault}{\updefault}$v_{(1)}$}}}
\put(1200,1300){\makebox(0,0)[lb]{{\SetFigFont{11}{14.4}{\rmdefault}{\mddefault}{\updefault}$u_{(1)}$}}}
\put(1400,3870){\makebox(0,0)[lb]{{\SetFigFont{11}{14.4}{\rmdefault}{\mddefault}{\updefault}$u_{(2)}$}}}
\put(2925,2350){\makebox(0,0)[lb]{{\SetFigFont{11}{14.4}{\rmdefault}{\mddefault}{\updefault}$u_{(1,2)}$}}}
\put(2970,1400){\makebox(0,0)[lb]{{\SetFigFont{11}{14.4}{\rmdefault}{\mddefault}{\updefault}$v_{(1,2)}$}}}
\put(-100,1500){\makebox(0,0)[lb]{{\SetFigFont{11}{14.4}{\rmdefault}{\mddefault}{\updefault}$v$}}}
\put(-100,2450){\makebox(0,0)[lb]{{\SetFigFont{11}{14.4}{\rmdefault}{\mddefault}{\updefault}$u$}}}
\put(500,500){\makebox(0,0)[lb]{{\SetFigFont{11}{14.4}{\rmdefault}{\mddefault}{\updefault}$x^i$}}}
\put(500,3250){\makebox(0,0)[lb]{{\SetFigFont{11}{14.4}{\rmdefault}{\mddefault}{\updefault}$y^i$}}}
\put(2400,3250){\makebox(0,0)[lb]{{\SetFigFont{11}{14.4}{\rmdefault}{\mddefault}{\updefault}$p^i$}}}
\put(2400,600){\makebox(0,0)[lb]{{\SetFigFont{11}{14.4}{\rmdefault}{\mddefault}{\updefault}$q^i$}}}
\end{picture}
} 
\caption{2-field lattice equations: the values of the fields are assigned on the vertices of a cube,
while the YB variables on the faces. }
\label{fig:orient}
\end{figure}
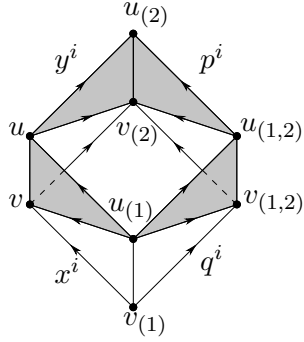

The above framework provides a method for obtaining a YB map from 
a three-dimensional consistent lattice system of the form (\ref{eq:quad}) $(\ell=2)$.
However, there are circumstances where invariance of a system of lattice equations under a 
one-parameter group of transformations suffices to construct an associated YB map $R$.
This is the case for all lattice equations of table \ref{tab1}, where we make the identification 
$f=(u,v)$. In order to visualize better this situation we are going to consider 
a cube with base the elementary square and place the second field on the top face 
as in figure \ref{fig:orient}. 

Let us examine first the system ${\mathcal E}_1$ which, for convenience, we rewrite below 
\begin{equation}
u_{(1,2)}=u+\frac{\beta_{1}-\beta_2 }{1-u_{(1)} v_{(2)}}\,u_{(1)} \,, \qquad
v_{(1,2)}=v-\frac{\beta_{1}-\beta_2 }{1-u_{(1)} v_{(2)}}\,v_{(2)} \,. \label{eq:AdlYam88}
\end{equation}
Notice that the values $u_{(2)}$, $v_{(1)}$ are not involved in the equations. 
Clearly, equations (\ref{eq:AdlYam88}) are invariant under the one-parameter group of transformations
\begin{equation}
G:\,\, (u,v) \mapsto (e^\varepsilon \, u ,e^{-\varepsilon} \, v)\,, \label{eq:scale1}
\end{equation}
with corresponding infinitesimal generator 
\begin{equation}
{\mathbf{v}}_1 = u\, \partial_u- \, v\, \partial_v \, .  \label{eq:vscale1}
\end{equation}
We consider now the space ${\rm J}^{(2)}(\stackrel{\vee}{v_{(1)}})$, i.e.
the space obtained from ${\rm J}^{(2)}$ by factoring out $v_{(1)}$. 
On ${\rm J}^{(2)}(\stackrel{\vee}{v_{(1)}})$ we
define the following lattice invariants along the orbits of ${\mathbf{v}}_1$
\begin{subequations} 
\label{eq:inveqs}
\begin{align}\label{eq:inveq1}
(x^1,x^2) = ({u}/{u_{(1)}}\, ,\,v\, u_{(1)})\,, & \qquad (p^1,p^2) = ({u_{(2)}}/{u_{(1,2)}}\, ,\,v_2\, u_{(1,2)})\, , \\
(y^1,y^2) = ({u_{(2)}}/{u}\, , \,u\, v_{(2)})\,, &  \qquad (q^1,q^2) =({u_{(1,2)}}/{u_{(1)}}\, , \,u_{(1)}\, v_{(1,2)})\,. \label{eq:inveq2}
\end{align}
\end{subequations} 
The space of $G$-orbits on ${\rm J}^{(2)}(\stackrel{\vee}{v_{(1)}})$ is six dimensional and 
thus there exist exactly two functionally independent relations
among the invariants (\ref{eq:inveqs}), namely
\begin{equation}
 p^1\, q^1 = x^1\, y^1\,, \quad p^1\, p^2 = y^1\, y^2 \,. \label{eq:functeq1}
\end{equation}
Moreover, since $G$ is a symmetry, equations (\ref{eq:AdlYam88}) can be written in terms of the 
invariants (\ref{eq:inveqs}), in the following form
\begin{equation}
 p^2 = y^2 P\,, \quad q^2=x^2 + y^2(1-P)\,,  \qquad \mbox{where} \qquad P = 1+ \frac{\beta_1-\beta_2}{x^1-y^2} \,.
 \label{eq:eq1invform}
\end{equation}
Solving the system of algebraic equations (\ref{eq:functeq1}), (\ref{eq:eq1invform}) 
for $(p,q)$ in terms of $(x,y)$, we obtain the unique solution
\begin{equation}
(p^1,p^2) = (y^1 \, P^{-1} \, , \, y^2 \, P )\,,  \qquad (q^1,q^2)=(x^1\, P \, , \, x^2 +y^2(1 - P ) \,) \,,
\label{eq:YBmapAY}
\end{equation}
which defines the YB map $R: \, (x^i,y^i) \mapsto (p^i,q^i)$ (shaded triangles of figure \ref{fig:orient}).

On the other hand, one could equally well consider on the space ${\rm J}^{(2)}(\stackrel{\vee}{u_{(2)}})$  
the following invariants under the action of $G$ (blank triangles of figure \ref{fig:orient})
\begin{subequations} 
\begin{align}
(x^3,x^2) = ({v}/{v_{(1)}}\, ,\,v\, u_{(1)})\,, & \qquad (p^3,p^2) = (v_{(2)}/v_{(1,2)}\, ,\,v_{(2)}\, u_{(1,2)})\,, \\
(y^3,y^2) = (v_{(2)}/v\, , \,u\, v_{(2)})\,, &  \qquad (q^3,q^2) = (v_{(1,2)}/v_1\, , \,u_{(1)}\, v_{(1,2)}) \,.
\end{align}
\end{subequations} 
The above invariants are related by 
\begin{equation}
 p^3\, q^3 = x^3\, y^3\,, \quad p^3\, q^2 = y^3\, x^2 \,, \label{eq:blankinv}
\end{equation}
while equations (\ref{eq:AdlYam88}) can be written in terms of them as follows
\begin{equation}
 p^2=y^2 + x^2\,(1 - \, \widetilde{P}) \,, \qquad q^2 =x^2\, \widetilde{P} \, , \qquad \mbox{where} \qquad 
\widetilde{P} = 1 - \frac{(\beta_1-\beta_2) y^3}{1- x^2\,y^3} \,.
\end{equation}
Solving the above equations for $(p^i,q^i)$, we obtain the YB map 
\begin{equation} (p^3,p^2) = (y^3 \, \widetilde{P}^{-1} \, , \, y^2 + x^2(1 - \, \widetilde{P}) )\,,  \qquad 
(q^3,q^2)=(x^3 \, \widetilde{P} \, , \, x^2\, \widetilde{P})\,. 
\end{equation}
The variables (\ref{eq:inveqs}) (here called Yang-Baxter variables) 
were introduced in a different context in \cite{A-Y}, where the reduction of the lattice equation (\ref{eq:AdlYam8}) to the 
map (\ref{eq:YBmapAY}) was also presented, together with an associated Lax pair. 
It seems that one of the motivations in \cite{A-Y} for this reduction  
was the illustration of the discrete Liouville theorem for integrable maps. 

We close this section by presenting a Lax pair, different from the one given 
in \cite{A-Y}, for the YB map (\ref{eq:YBmapAY}). The YB map (\ref{eq:YBmapAY})
admits the discrete zero curvature representation
\begin{equation}
W(y^1,y^2,\beta_2,\lambda)\, W(x^1,x^2,\beta_1,\lambda)\,=
W(p^1,p^2,\beta_1,\lambda)\,W(q^1,q^2,\beta_2,\lambda)\,, \nonumber
\end{equation}
where 
\begin{equation}W(\xi^1,\xi^2,\beta,\lambda)=\left(
\begin{array}{cc}\xi^1 & - \xi^1 \, \xi^2 \\ 1 & \lambda -\beta -\xi^2 \end{array}
\right) \,. \nonumber
\end{equation}
A hint for considering the matrix $W$ of the above form, is based on the YB map itself. Indeed,
in view of equations (\ref{eq:YBmapAY}) giving $p^1$, one notices that
\begin{equation}
p^1 = y^1 \, \frac{x^1-y^2}{x^1-y^2+\beta_1-\beta_2}\,, \nonumber
\end{equation}
can be written as a linear fractional transformation induced from the linear transformation
with matrix $W(y^1,y^2,\beta_2,\beta_1)$ on $(x^1,1)$ (column vector), cf \cite{Suris}.

\section{YB maps arising from equations $\mathcal{E}_{\alpha}$}
\label{sec5}

In this section we present the YB maps obtained from the lattice equations ${\mathcal{E}}_\alpha$ by
considering invariants of their symmetry groups. For equations admitting a two-parameter
symmetry group we consider as YB variables a complete set of joint invariants on $\rm{J}^{(2)}$.
In the case of a one-parameter symmetry group, as we discussed in the previous section, one can consider
two cases for YB variables: (i) lattice invariants on the space ${\rm J}^{(2)}(\stackrel{\vee}{v_{(1)}})$ 
(shaded triangles in figure \ref{fig:orient}), or (ii) lattice invariants on the space ${\rm J}^{(2)}(\stackrel{\vee}{u_{(2)}})$ 
(blank triangles in figure \ref{fig:orient}). We present YB maps obtained from case (i) only. 
For each system the analysis of case (ii) yields identical map as in case (i),
 up to appropriate rearrangement of variables.

\subsection{System ${\mathcal{E}}_2$}

Equations $\mathcal{E}_2$ admit the symmetry generator (\ref{eq:vscale1}), thus we take as YB 
variables the ones defined by equations (\ref{eq:inveqs}), leading to the same functional
relations (\ref{eq:functeq1}). On the other hand, system $\mathcal{E}_2$ can be written in terms of
the invariants (\ref{eq:inveqs}) as 
\begin{equation}
p^2 = y^2 \, P \, ,\qquad  q^2 = x^2 \, Q\,, \label{eq:G2rel}
\end{equation}
where 
\begin{equation}
\displaystyle{P = 1+ (\beta_2-\beta_1) \frac{1-x^1}{\beta_2 \, x^1 + y^2} }\,, \qquad
\displaystyle{Q=1+(\beta_1-\beta_2) \frac{y^2/x^2-x^1\, }{\beta_1 \, x^1 + y^2}}\,.  \nonumber
\end{equation}
Solving the system of equations (\ref{eq:functeq1}), (\ref{eq:G2rel}) for $(p^i,q^i)$ 
we obtain the YB map
\begin{equation}
(p^1,p^2) = (y^1 \, P^{-1}\,, y^2 \, P) \,, \quad (q^1,q^2)= (x^1\, P \,, x^2 \, Q )\,. \nonumber
\end{equation}
\subsection{System ${\mathcal{E}}_3$}
 Equations $\mathcal{E}_3$ admit the symmetry generator
\begin{equation}
 \mathbf{v}_2=\partial_u - \partial_v\,. \label{eq:translgen}
\end{equation}
On the space ${\rm J}^{(2)}(\stackrel{\vee}{v_{(1)}})$ we define as YB variables the following invariants
\begin{subequations} 
\label{eq:inv1G3all}
\begin{align}
(x^1,x^2) = (u-u_{(1)} \, ,\,v + u_{(1)})\,,  \quad (p^1,p^2) = (u_{(2)} - u_{(1,2)}\, ,\,v_{(2)} + u_{(1,2)})\,, \label{eq:inv1G3}  
\\ (y^1,y^2) = (u_{(2)} - u\, , \,u + v_{(2)})\,,  
 \quad (q^1,q^2) = (u_{(1,2)} - u_{(1)}\, , \,u_{(1)} +  v_{(1,2)})\,.
\label{eq:inv2G3} 
\end{align}
\end{subequations} 
They are related by
\begin{equation}
p^1 +  q^1 = x^1 +  y^1\,, \qquad 
p^1 + p^2 = y^1 + y^2\,,  \label{eq:G3eq1}
\end{equation}
and system $\mathcal{E}_3$ is written in terms of the invariants (\ref{eq:inv1G3all}) as follows
\begin{equation}
p^2 = y^2 \, P\,, \qquad q^2 = x^2 \, P^{-1} \,, \label{eq:G3eq2}
\end{equation}
where 
\begin{equation} 
\displaystyle{P = \frac{y^2 - x^1 -\beta_2}{y^2 - x^1 -\beta_1} }\,.  \nonumber
\end{equation}
Solving the system of equations (\ref{eq:G3eq1}), (\ref{eq:G3eq2}) for $(p^i,q^i)$ we obtain the YB map
\begin{equation} 
(p^1,p^2) = \big(y^1 + y^2 \, ( 1 - \, P) \,, \, y^2 \, P\big)\,, \quad
(q^1,q^2) = \big(x^1-y^2 \,(1 -P)\,, \, x^2 \,  P^{-1}\big)\,. \nonumber 
\end{equation}
\subsection{System ${\mathcal{E}}_4$}
In the generic case $\epsilon \neq 0$, the YB variables are the invariants 
of the symmetry group generated by the vector field (\ref{eq:translgen}). System $\mathcal{E}_4$
is written in terms of the invariants (\ref{eq:inv1G3all}) as follows
\begin{equation}
p^2 = y^2 - P\, ,\, \quad q^2 = x^2 + Q\,,  \label{eq:G4eq2}
\end{equation}
where 
\begin{equation} 
\displaystyle{P= (\frac{1}{\beta_1}-\frac{1}{\beta_2}) \, \frac{\beta_1\,\beta_2\,y^2 \, x^1 + {\epsilon}}
{\beta_2\, y^2 - \beta_1 \, x^1} }\,, \quad
\displaystyle{Q=  (\frac{1}{\beta_2}-\frac{1}{\beta_1}) \, \frac{ \beta_1\,\beta_2\, x^2\,(y^2-x^1-x^2) - 
{\epsilon}}{\beta_1\, x^2 + \beta_2 \, (y^2-x^1-x^2)} } \,.  \nonumber
\end{equation}
Inserting equations (\ref{eq:G4eq2}) into (\ref{eq:G3eq1}) and solving for $(p^1,q^1)$ in terms of $(x^i,y^i)$, 
we arrive at the YB map 
\begin{equation}
(p^1,p^2) = (y^1 + P\, ,\, y^2 - P )\,, \qquad 
(q^1,q^2) = (x^1 - P\, ,\,  x^2 + Q) \,. \nonumber
\end{equation}
We examine now the case $\epsilon=0$. In this case there exists, in addition, the scaling symmetry 
generated by the vector field
\begin{equation}
\mathbf{v}_3=u\,\partial_u + v\,\partial_v\,. \label{eq:scalG4}
\end{equation}
On ${\rm J}^{(2)}(\stackrel{\vee}{v_{(1)}})$ we define as YB variables the following invariants under $\mathbf{v}_3$
\begin{subequations} 
\label{eq:invG40all}
\begin{align}
{ (x^1,x^2) = ({u}/{u_{(1)}} \, ,\, {v}/{u_{(1)}} ) }\,,  \quad
{ (p^1,p^2) = ({u_{(2)}}/{u_{(1,2)}}\, ,\,{v_{(2)}}/{u_{(1,2)}}) }\,, \label{eq:invG40eq1}  \\
{ (y^1,y^2) = ({u_{(2)}}/{u}\, , \, {v_{(2)}}/{u}) }\,, \quad
{ (q^1,q^2) = ({u_{(1,2)}}/{u_{(1)}} \, , \, {v_{(1,2)}}/{u_{(1)}} ) } \,, \label{eq:invG40eq2}
\end{align}
\end{subequations}
which are related by
\begin{equation} 
p^1 \, q^1 = x^1 \, y^1 \, , \quad y^2 \, p^1 = y^1 \,  p^2 \,. \label{eq:funcG40}
\end{equation}
System ${\mathcal{E}}_4$ ($\epsilon=0$) can be written in terms of the invariants (\ref{eq:invG40all}), as follows
\begin{equation} 
p^1 = x^1 \, P\, ,\, \qquad q^2 = x^2 \, Q\,,  \label{eq:G40invform}
\end{equation}
where
\begin{eqnarray}
\displaystyle{ P =  1+  \frac{(\beta_2-\beta_1) \,(1+y^2)\,(1-x^1)}{\beta_2\, x^1 \, (1+y^2) + \beta_1 \, (1-x^1)} }\,,  \nonumber \\
\displaystyle{ Q= 1+ \frac{1}{x^2} \, \frac{(\beta_1-\beta_2)\,(1+x^2)\,(x^1 \, y^2 - x^2)}{\beta_1\, (1+x^2) + 
\beta_2 \, ( x^1 \, y^2-x^2)} }\,.  \nonumber
\end{eqnarray}
From relations (\ref{eq:funcG40}), (\ref{eq:G40invform}) we obtain the following YB map
\begin{equation} 
(p^1,p^2) = (y^1 \, P^{-1}\, , \, y^2 \, P^{-1}) \,, \qquad 
(q^1,q^2) = (x^1 \, P \,, \, x^2 \, Q ) \,. \nonumber
\end{equation}
Finally, considering the group generated by 
$\lbrace \mathbf{v}_2, \mathbf{v}_3 \rbrace$, given by
(\ref{eq:translgen}), (\ref{eq:scalG4}), we take on ${\rm{J}}^{(2)}$ 
the following joint invariants
\begin{subequations}
\label{eq:jointall}
\begin{align} 
{ (x^1,x^2) = \big(\frac{u-u_{(1)}}{v+u_{(1)}} \, ,\,\frac{v-v_{(1)}}{v+u_{(1)}} \big) }  ,\quad
{ (p^1,p^2) = \big(\frac{u_{(2)}-u_{(1,2)}}{v_{(2)}+u_{(1,2)}}\, ,\,\frac{v_{(2)}-v_{(1,2)}}{v_{(2)}+u_{(1,2)}}\big) } \, , 
\label{eq:joint1} \\  
{ (y^1,y^2) = \big(\frac{u_{(2)}-u}{u+v_{(2)}}\, , \, \frac{v_{(2)}-v}{u+v_{(2)}}\big) }\, , \quad
{ (q^1,q^2) = \big(\frac{u_{(1,2)}-u_{(1)}}{u_{(1)}-v_{(1,2)}} \, , \, \frac{v_{(1,2)}-v_{(1)}}{u_{(1)}+v_{(1,2)}} \big) } \,. \label{eq:joint2}
\end{align}
\end{subequations}
The above invariants are functionally related by
\begin{equation}
\displaystyle{\frac{(1+ x^1) \, (1+y^1)}{(1+p^1)\,(1+q^1)} = \frac{(1-x^2)\, (1-y^2)}{ (1-p^2) \, (1-q^2) } =
\frac{1+ x^1 \, y^2 }{1+ p^2\, q^1} } \,. \label{eq:joint1rel}
\end{equation}
In terms of the invariants (\ref{eq:jointall}), the system $\mathcal{E}_4$ ($\epsilon=0$) 
can be written  as follows
\begin{equation}
q^1=x^1 \, P\, ,  \quad p^2=y^2 \, P^{-1}\,, \label{eq:jointG40}
\end{equation}
where
\begin{equation}
\displaystyle{P= \frac{ \beta_1\,(1-y^2)+ \beta_2\, y^2\, (1+x^1)}{\beta_1\,x^1\,(y^2-1)+\beta_2\,(1+x^1)}} \,. \nonumber
\end{equation}
Using relations (\ref{eq:joint1rel}), (\ref{eq:jointG40}) we obtain the YB map
\begin{subequations}
\label{eq:YBmap40all}
\begin{align}
&\displaystyle{ (p^1,p^2) = \left( \frac{(1+x^1)\,(1+y^1)}{1+ x^1 \, P} -1  \, , \,  y^2 \,  P^{-1} \right) \,, } \label{eq:YBmap401} \\
&\displaystyle{ (q^1,q^2) = \left( x^1 \, P \,, \, \frac{(1-x^2)\,(1-y^2)} {y^2\,P^{-1} - 1} +1  \right) }\,. \label{eq:YBmap402}
\end{align}
\end{subequations}
\subsection{System ${\mathcal{E}_5}$} 
We examine first the generic case $\delta=1$, i.e.
 $\gamma_i^2-\beta_i^2=1$. System $\mathcal{E}_5$ admits
the scaling symmetry with generator $\mathbf{v}_3$, thus we use the invariants (\ref{eq:invG40all}). 
In terms of these variables the system $\mathcal{E}_5$ is written as follows
\begin{equation}
p^1 = x^1 \, P\, ,\, \qquad q^2 = x^2 \, Q\,, \label{eq:G5inv}
\end{equation}
where
\begin{equation}
\displaystyle{ P =  \frac{(\beta_1-\beta_2) \,y^2 + \gamma_2 + \gamma_1 \, x^1 \, y^2}{(\beta_1-\beta_2) \, x^1 + \gamma_1 + \gamma_2 \, x^1 \, y^2} } \,, \quad
\displaystyle{ Q= \frac{(\beta_2-\beta_1) \, \frac{x^1\,y^2}{x^2} + \gamma_2 + \gamma_1 \, x^1\, y^2}{(\beta_2-\beta_1)\, x^2 + \gamma_1 + \gamma_2 \,  x^1 \, y^2} } \,. \nonumber
\end{equation}
The YB map obtained from relations (\ref{eq:funcG40}), (\ref{eq:G5inv}) is
\begin{equation}
(p^1,p^2) = (y^1 \,  P^{-1} \,,\, y^2 \, P^{-1} ) \,, \qquad
(q^1,q^2) = (x^1 \, P \, , \, x^2 \, Q)\,.  \nonumber
\end{equation}
When $\delta=0$, the system $\mathcal{E}_5$ becomes system
$\mathcal{E}_4$ with $\epsilon=0$ directly (case $\gamma_i=-\beta_i$), or
indirectly under the gauge transformation $(u,v)\mapsto (-1)^{n_1+n_2}(u,v)$ 
(case $\gamma_i=\beta_i$). 

\section{Discrete vector Schr\"odinger system}
\label{sec6}
In this section we show how the preceding considerations can be generalized to vector extensions
of integrable lattice equations yielding $n$--vector (two-field) YB maps. This is illustrated for the lattice
$n$-vector Schr\"odinger system 
\begin{equation}
\boldsymbol{u}_{(1,2)} = \boldsymbol{u} + \frac{\beta_1-\beta_2}{1-\boldsymbol{u}_{(1)} \cdot \boldsymbol{v}_{(2)}}\,\boldsymbol{u}_{(1)} \,, \qquad
\boldsymbol{v}_{(1,2)} = \boldsymbol{v} + \frac{\beta_2-\beta_1}{1-\boldsymbol{u}_{(1)} \cdot \boldsymbol{v}_{(2)}}\,\boldsymbol{v}_{(2)}\,,
\label{eq:vNLS}
\end{equation}
introduced by Adler \cite{Adler}, where $\boldsymbol{u}=(u^1,\ldots,u^n)$ etc and $\cdot$ stands for
the usual Euclidean inner product. 
Equations (\ref{eq:vNLS}) represent the superposition formula of the B\"acklund transformation for the
$n$--vector nonlinear Schr\"odinger system found by Manakov \cite{Manakov}.

System (\ref{eq:vNLS}) remains invariant under the group of transformations generated by the vector fields 
\begin{equation}
\mathbf{v}_i = {u^i} \, \partial_{u^i} - {v^i} \, \partial_{v^i}\,, \nonumber
\end{equation}
$i=1,\ldots n$, where no summation on indices is assumed.
In order to make the presentation more concise, we represent the vectors $\boldsymbol{u}, \boldsymbol{v}$ 
and their shifts in terms of diagonal matrices denoted by the corresponding capital letter, e.g. 
$\boldsymbol{u} \rightarrow  U ={\rm{diag}}\,(u^1,u^2,\ldots, u^n) $. On the lattice jet space 
${\rm J}^{(2)}(\stackrel{\vee}{\boldsymbol{v}_{(1)}})$ we consider the following YB variables
\begin{align}
&\big(X^1,X^2\big) = \big(\,{U} \, {U_{(1)}}^{-1}\, ,\,V\, {U_{(1)}}\,\big)\,,  \qquad
\big(P^1,P^2\big) = \big(\,{U_{(2)}}\,{U_{(1,2)}}^{-1}\, ,\,{V_{(2)}}\, {U_{(1,2)}}\,\big)\, , \nonumber  \\
&\big(Y^1,Y^2\big) = \big(\,{U_{(2)}}\,U^{-1}\, , \,U\, {V_{(2)}}\,\big)\, , \qquad\,\,
\big(Q^1,Q^2\big) = \big(\,{U_{(1,2)}}\,{U_{(1)}}^{-1}\, , \,{U_{(1)}}\, {V_{(1,2)}}\,\big)\, . \nonumber
\end{align}
The above invariants are functionally related by
\begin{equation} 
P^1 \, Q^1 = X^1 \, Y^1 \,,\qquad P^1 \, P^2 = Y^1 \, Y^2 \,, \label{eq:vNLSrel1}
\end{equation}
and the lattice equations (\ref{eq:vNLS}) are written in terms of them as follows
\begin{equation}
P^2 = Y^2 \, S\,,  \qquad  P^2 + Q^2  = X^2 + Y^2 \,, \label{eq:vNLSrel2}
\end{equation}
where 
\begin{equation}
S = I + \frac{(\beta_1-\beta_2)\, {{(X^1)}^{-1}}}{1- {\rm tr} \big(Y^2\, {(X^1)}^{-1}\big)} \,. \nonumber
\end{equation}
Solving equations (\ref{eq:vNLSrel1}), (\ref{eq:vNLSrel2}) for $(P^i,Q^i)$ in terms of $(X^i,Y^i)$ we
obtain the YB map 
\begin{equation} \big(P^1,P^2\big) = \big(Y^1 \, S^{-1}\, ,\, Y^2\,S\,\big) \,,
\quad \big(Q^1,Q^2\big) = \big( X^1 \, S \, , \, X^2 + Y^2 \, (I-S)  \,\big)\,. \label{eq:YBvNLS}
\end{equation}
It would be interesting to investigate the relation of the YB maps studied recently 
in \cite{Abl1,Tsuchida,Abl2} with the YB map (\ref{eq:YBvNLS}) after imposing
appropriate reality conditions on the latter, as it was done for system 
${\mathcal{E}}_1$ in \cite{A-Y}.

\section{Conclusions and perspectives}
\label{persp}
We have applied symmetry methods to multi-field integrable lattice
equations in order to construct YB maps. The main idea is that the invariants of the symmetry
groups of the lattice equations serve as variables of the associated YB maps.

According to \cite{ABS1}, a birational map $ R: (x,y) \mapsto  (p(x,y),q(x,y))$ is called 
quadrirational if the maps $s \mapsto p(s,y)$ and $t \mapsto q(x,t)$ for generic values of $y$, $x$, 
are invertible. This property is equivalent to the nondegeneracy property which is often imposed 
additionally in the study of YB maps, cf \cite{ESS}.
We have to notice that all YB maps we presented in this work do not satisfy 
the quadrirationality property. This can be seen immediately since $(p^1,p^2)$ is 
independent of $x^2$. Two additional properties which characterize certain YB maps
are invertibility and unitarity. A YB map $R(\beta_1,\beta_2)$ with parameters $(\beta_1,\beta_2)$
is called unitary if it satisfies the relation 
\begin{equation}
R^{(2,1)}(\beta_2,\beta_1)\,R(\beta_1,\beta_2)=\rm{id}\,, \nonumber
\end{equation}
where $R^{(2,1)} = P\, R \, P$ and $P$ the permutation map i.e. $P:(x,y)\mapsto (y,x)$. 
The YB maps we obtained here are invertible and unitary.

Recently, Reshetikhin and Veselov \cite{NS} studied the Hamiltonian theory 
of YB maps through Lie-Poisson structures. 
The Poisson structure corresponding to the auto-transformation
(\ref{eq:AdlYam8}) was presented in \cite{A-Y}. It would be interesting to explore
the Poisson structure of the remaining auto-transformations and the associated
YB maps using the theory of Poisson Lie groups as in \cite{WX,NS}.

The protagonistic role of the YB relation in quantum integrability is well-known. 
In view of recent developments (see \cite{ves2} and references therein) 
it seems that the YB relation for maps acquires a 
prominent role in (classical) discrete integrability as well. 
Here, we saw that integrable lattice equations admitting a symmetry group
can be reduced to YB maps.
It remains an open question how integrable lattice equations which do not admit an
appropriate symmetry group can be casted in the form of a YB map. 

By the way of studying the multi-field lattice equations introduced by Adler and Yamilov the
question of classification of two-field lattice equations by imposing the three-dimensional
consistency property is naturally raised. Before tackling this problem in full generality
one can consider the restricted class of equations of type (\ref{eq:AdlYamR}). 
It would be interesting to investigate whether the 
classification is exhausted essentially by the lattice equations in \cite{A-Y}.

\subsection*{Acknowledgements}
VGP acknowledges support from the programme ``C.
Carath\'eodory" of the Research Committee of the University of
Patras. The work of AGT was supported by the 
grant Pythagoras B-365-015 of the European Social Fund  (EPEAEK II).

\end{document}